\documentclass[a4paper,11pt]{amsart}
\usepackage{wrapfig}
\usepackage[dvips]{graphicx}
\usepackage{amsmath,amsthm,amssymb,amscd}
\theoremstyle{definition}
\newtheorem{thm}{Theorem}[section]

\newtheorem{pro}[thm]{Proposition}
\newtheorem{cor}[thm]{Corollary}
\newtheorem{lem}[thm]{Lemma}
\newtheorem{ex}[thm]{Example}
\newtheorem{rem}[thm]{Remark}
\newtheorem{fac}[thm]{Fact}
\title{$C^*$-algebras associated with real multiplication}
\author{Norio Nawata}
\address{Graduate~School~of~Mathematics,~Kyushu~University,~Hakozaki,~Fukuoka,~812-8581,~Japan}
\email{n-nawata@math.kyushu-u.ac.jp}
\begin{document}
\begin{abstract}
Noncommutative tori with real multiplication are the irrational rotation algebras that have special
equivalence bimodules. Y. Manin proposed the use of noncommutative tori with real multiplication
as a geometric framework for the study of abelian class field theory of real quadratic fields.
In this paper, we consider the Cuntz-Pimsner algebras constructed by special equivalence bimodules 
of irrational rotation algebras.  We shall show that associated $C^*$-algebras are simple and
purely infinite. We compute the K-groups of associated $C^*$-algebras and show that these algebras are related
to the solutions of Pell's equation and the unit groups of real quadratic fields. We consider the 
Morita equivalent classes of associated $C^*$-algebras.
\end{abstract}
\maketitle
\section{Introduction}
Let $\theta$ be an irrational number. An irrational rotation algebra $A_\theta$ is the 
crossed product $C^*$-algebra for the action of the integers on the circle
by powers of the rotation by angle $2\pi\theta$. It is simple and has a unique normalized trace 
$\tau_\theta$. These algebras are also called noncommutative tori and have been classified up to 
$C^*$-isomorphism and Morita equivalence \cite{Pim1},\cite{Rie1}. 

Suppose that $\sigma$ be a free and proper action of a locally compact group $G$ on a 
locally compact Haussdorff space $X$. Then $X/\sigma$ is a locally compact Haussdorff space 
and $C_0(X)\rtimes_\sigma G$ is Morita equivalent to $C_0(X/\sigma) $ (see \cite{RieG}). 
The action of the integers on the circle by powers of the rotation by angle $2\pi\theta$ 
is free, but every orbit is dense and the orbit space $(\mathbb{R}/\mathbb{Z})/\mathbb{Z}\theta$ 
behaves badly. 
By the fact above, we may consider an irrational rotation algebra $A_\theta$ serves to 
represent a quotient space $(\mathbb{R}/\mathbb{Z})/\mathbb{Z}\theta$. 
This idea is due to A. Connes \cite{Con}. 

An elliptic curve $E_\omega$ over $\mathbb{C}$ can be described as the quotient 
$E_\omega\cong\mathbb{C}/(\mathbb{Z}+\mathbb{Z}\omega)$ 
of the complex plane by a 2-dimensional lattice $\mathbb{Z}+\mathbb{Z}\omega$ where we can 
take $Im(\omega )>0$. Most elliptic curves over $\mathbb{C}$ have only the multiplication-by-$m$ 
endomorphisms. An elliptic curve $E_\omega$ over $\mathbb{C}$ has extra endomorphisms if and 
only if $\omega$ is in a imaginary quadratic field. In this case, an elliptic curve $E_\omega$ 
is said to have complex multiplication. Such curves have many special properties. One of the 
important properties is that $j$-invariant and torsion points of $E_\omega$ generate a maximal 
abelian extension of a imaginary quadratic field $\mathbb{Q}(\omega)$. For real quadratic field, 
a similar description is not known. By the idea above and regarding 
$(\mathbb{R}/\mathbb{Z})/\mathbb{Z}\theta$ as $\mathbb{R}/(\mathbb{Z}+\mathbb{Z}\theta)$, 
Y. Manin proposed to use irrational rotation algebras corresponding to real quadratic 
irrationalities as a replacement of elliptic curves with complex multiplication. 
We shall explain there exist the analogous properties with complex multiplication in 
the irrational rotation algebras.

Two $C^*$-algebras $A$ and $B$ are Morita equivalent if there exists an $A$-$B$-equivalence 
bimodule. An $A$-$B$-equivalence bimodule is an $A$-$B$-bimodule $\mathcal{E}$ which is 
simultaneously a full left
Hilbert $A$-module under an $A$-valued inner product $_A\langle\cdot ,\cdot\rangle$ and a full
right Hilbert $B$-module under a $B$-valued inner product $\langle\cdot ,\cdot\rangle_B$, satisfying
$_A\langle\xi ,\eta\rangle\zeta =\xi\langle\eta ,\zeta\rangle_B$ for any 
$\xi ,\eta , \zeta\in\mathcal{E}$. Let $\alpha$ be an 
automorphism of a $C^*$-algebra $A$. We denote by $\mathcal{E}_{\alpha}$ the vector space $A$
with the obvious left $A$-action and the obvious left $A$-valued inner product. We define the 
right
$A$-action on $\mathcal{E}_\alpha$ by $\xi\cdot a=\xi\alpha^{-1}(a)$ for any $\xi\in
\mathcal{E}_\alpha$ and $a\in A$, and the right $A$-valued inner product by $\langle\xi ,\eta\rangle_A
=\alpha (\xi^*\eta)$ for any $\xi ,\eta\in\mathcal{E}_\alpha$. Then $\mathcal{E}_\alpha$ is an
$A$-$A$-equivalence bimodule. Hence equivalence bimodules are regarded as a generalization 
of automorphisms of $C^*$-algebras. The Picard group of $C^*$-algebra $A$ is the set of isomorphic classes 
of the $A$-$A$-equivalence bimodules \cite{Bro}. It forms a group structure whose product is defined by tensor 
product. Let $\alpha$ and $\beta$ be automorphisms of $A$. Then $\mathcal{E}_\alpha$ is isomorphic
to $\mathcal{E}_\beta$ if and only if $\alpha$ is unitary equivalent to $\beta$. Moreover, 
$\mathcal{E}_\alpha \otimes \mathcal{E}_\beta$ is isomorphic to $\mathcal{E}_{\beta\circ\alpha}$.
Hence we obtain an anti-homomorphism of $\mathrm{Aut} (A)/\mathrm{Inn} (A) $ to the Picard group 
of $A$.   

The equivalence bimodules of irrational rotation algebras are constructed by M. Rieffel. 
These bimodules are defined as completions of $C_c(\Bbb{R}\times\Bbb{Z} /c\Bbb{Z} )$ or $S
(\Bbb{R}\times\Bbb{Z} /c\Bbb{Z} )$ with the certain actions of irrational rotation algebras 
\cite{Rie2}, \cite{Rie3}. There exists 
special $A_\theta$-$A_\theta$-equivalence bimodules in the case where $\theta$ is a quadratic irrational number. 
In fact, K. Kodaka showed that if 
$\theta$ is not a quadratic irrational number, then the Picard group of $A_\theta$ is isomorphic
to $\mathrm{Aut}(A_\theta )/\mathrm{Inn}(A_\theta )$ and that if $\theta$ is a quadratic number, then the Picard
group of $A_\theta$ is isomorphic to a semidirect product of $\mathrm{Aut} (A_\theta )/\mathrm{Inn} (A_\theta )$
with $\Bbb{Z}$ \cite{kod2}. We may consider that this is analogous to complex multiplication 
of elliptic curve. Y. Manin called this real multiplication and proposed the use of 
noncommutative tori with real multiplication as a geometric framework for the study of abelian 
class field theory of real quadratic fields \cite{Man}.

In this paper, we consider the Cuntz-Pimsner algebras constructed by the equivalence bimodules 
of irrational rotation algebras, that are not generated by automorphisms and show that these algebras 
are related to the solutions of Pell's equation and the unit groups of real quadratic fields. 

In Section \ref{sec:Rea} we consider the $A_\theta$-$A_\theta$-equivalence bimodules. 
K. Kodaka consider the $A_\theta$-$A_\theta$-equivalence bimodules as the automorphisms of 
$A_\theta\otimes\mathbb{K}$ where $\mathbb{K}$ is the $C^*$-algebras of all compact operators 
on a countably infinite dimensional Hilbert space and determine the Picard groups of the irrational 
rotation algebras in \cite{kod2}. We consider the equivalence bimodules of irrational 
rotation algebras neither as the completions of $C_c(\Bbb{R}\times\Bbb{Z} /c\Bbb{Z} )$ 
($S(\Bbb{R}\times\Bbb{Z} /c\Bbb{Z} )$) nor as the automorphisms of $A_\theta\otimes\mathbb{K}$ 
in this paper. 
If $A$ is unital, then an $A$-$A$-equivalence bimodule is a finitely generated 
projective $A$-module as a right $A$-module. Note that we need to be careful with a left action. 
We consider the $A_\theta$-$A_\theta$-equivalence 
bimodules from this viewpoint. But we analyze the $A_\theta$-$A_\theta$-equivalence bimodules by 
the similar arguments in \cite{kod2}. 

In Section \ref{sec:Ass} we construct $C^*$-algebras from the bimodules considered in Section \ref{sec:Rea}. 
We show that these algebras are purely infinite, simple, nuclear and in the UCT class.

In Section \ref{sec:K} we compute the K-groups of the associated $C^*$-algebras and show that 
these groups are related to the solutions of Pell's equation, where Pell's equation is $x^2-Dy^2=\pm 4$ 
for given an integer $D$, to be solved in integers $x$ and $y$. Moreover, we consider the 
Morita equivalent classes of associated $C^*$-algebras.

\section{Real multiplication}\label{sec:Rea}
In this section, we shall consider $A_\theta$-$A_\theta$-equivalence bimodules and review 
the result in \cite{kod2} by the slightly different viewpoint. 
We refer the reader
to \cite{Bla} for basic facts of $C^*$-algebras and equivalence bimodules. 

We shall review equivalence bimodules of unital $C^*$-algebras. 
Let $A$ be a unital $C^*$-algebra. Then an $A$-$A$-equivalence bimodule $\mathcal{E}$
is a finitely generated projective $A$-module as a right module. Hence $\mathcal{E}$ is isomorphic to
$qA^n$ as a right Hilbert $A$-module, where $q$ is a projection in $M_n(A)$. A right Hilbert 
$A$-module $qA^n$ has a structure of $qM_n(A)q$-$A$-equivalence bimodule with the obvious
left $qM_n(A)q$-action and the $qM_n(A)q$-valued inner product such that 
$_{qM_n(A)q}\langle q(a_i)_i,q(b_i)_i\rangle =q(a_jb_k^*)_{j,k}q$
for any $a_i,b_i\in A$. Since $\mathcal{E}$ is an $A$-$A$-equivalence bimodule, $qM_n(A)q$ is 
isomorphic to $A$. Note that the left $A$-action on $qA^n$ and the left $A$-valued inner 
product is dependent on an isomorphism $\phi$ from $A$ to $qM_n(A)q$.
Let $\alpha$ be an outer automorphism of $A$ and $\mathcal{F}$ be $qA^n$ with
the left $A$-action such that $a\cdot\xi=\phi\circ\alpha (a)\xi$ for any $a\in A$ and $\xi\in qA^n$. 
Then $\mathcal{E}$ is not isomorphic to $\mathcal{F}$ as an $A$-$A$-equivalence
bimodule. Therefore we need to be careful with an isomorphism from $A$ to $qM_n(A)q$.

We shall consider the condition of $q$ that $qA_\theta^n$ is an $A_\theta$-$A_\theta$-equivalence bimodule. 
We need the following well-known fact.
\begin{lem}\label{lem:com}
Let $q_1$, $q_2$ be projections in $A_\theta$. If $\tau_\theta(q_1)\geq\tau_\theta(q_2)$ 
(resp. $\tau_\theta(q_1)=\tau_\theta(q_2)$), then there exists a unitary element $w$ in 
$A_\theta$ such that $q_1\geq w^*q_2w$ (resp. $q_1= w^*q_2w$).
\end{lem}
Let $Tr_\theta :=\tau_\theta\otimes Tr$ be the unnormalized trace on 
$M_n(A_\theta )$ where $Tr$ is the usual trace on $M_N(\Bbb{C})$. We have the following lemma (see \cite{kod}).
\begin{lem}\label{lem:kod}
If $q$ is a proper projection in $M_n(A_\theta )$ such that $Tr_\theta(q)=k(c\theta+d)$ where
 $k$ is a natural number and $c$, $d$ are integers such that $gcd(c,d)=1$, then
\[qM_n(A_\theta )q\cong M_k(A_{\frac{a\theta +b}{c\theta +d}}) \]
for any $a, b\in \Bbb{Z}$ such that $ad-bc=\pm1$.
\end{lem}
Since $qA_\theta^n$ is a $qM_n(A_\theta )q$-$A_\theta$-equivalence bimodule, the condition of 
$q$ is equivalent to the condition that 
$Tr_\theta (q)=c\theta +d$ where there exist integers $a$ and $b$ such that 
$\frac{a\theta +b}{c\theta +d}=\theta$ by the lemma above. 
We consider the case where $\theta$ is not a quadratic number. The following proposition is Corollary 9 in \cite{kod2}. 
\begin{pro}
Let $\theta$ be an irrational number and $\mathcal{E}$ be an $A_\theta$-$A_\theta$-equivalence 
bimodule. Assume that $\theta$ is not a quadratic number. Then there exists an automorphism $\alpha$
such that $\mathcal{E}$ is isomorphic to $\mathcal{E}_\alpha$.
\end{pro}
\begin{proof}
Let $n$ be a natural number and $q$ be a projection in $M_n(A_\theta )$ such that $\mathcal{E}$ 
is isomorphic to $qA_\theta^n$. We denote by $\phi$ an isomorphism from $A_\theta$ to $qM_n(A_\theta )q$  
such that the left $A_\theta$-action on $qA_\theta^n$ is defined by $a\cdot \xi =\phi (a)\xi$ 
for any $a\in A_\theta$ and $\xi\in qA_\theta^n$. 
By the discussion above, there exist integers $a,b,c,d$ such that $Tr_\theta (q)=c\theta +d$ 
and $\frac{a\theta +b}{c\theta +d}=\theta$. 
Since $\theta$ is not a quadratic number, $Tr_\theta (q)=c\theta +d=1$. Hence, there exists a 
unitary element $w$ in $M_n(A_\theta )$ such that $w^*qw=1\otimes e_{11}$ where $e_{11}$ is a 
rank one projection in $M_n(\Bbb{C})$ by Lemma \ref{lem:com}. Since $w^*qwA_\theta^n$ where 
the left $A_\theta$-action is defined by $a\cdot \xi =w^*\phi (a)w\xi$ for any $a\in A_\theta$ and 
$\xi\in w^*qwA_\theta^n$ is isomorphic to $qA_\theta^n$, $\mathcal{E}$ is isomorphic to 
$(1\otimes e_{11})A_\theta^n$. It is easy to see that $(1\otimes e_{11})A_\theta^n$ is 
isomorphic to $A_\theta$ where the left $A_\theta$-action is defined by $a\cdot b=w^*\phi (a)wb$ for any 
$a,b\in A_\theta$. Therefore the proof is complete.
\end{proof}
We consider the case where $\theta$ is a quadratic irrational number. We may assume that 
$\theta$ satisfies $k\theta^2+l\theta+m=0$ with a natural number $k$ and integers $l$,$m$ 
such that $gcd(k,l,m)=1$. The equation is uniquely determined. Let $D_\theta =l^2-4km$ be the 
discriminant of $\theta$. The following fact is well-known (see \cite{gau} and \cite{kod}(Remark 7)).
\begin{fac}
Let $\theta$ be a quadratic number. Then there exists a real number $\epsilon_0$ such that
\[\{c\theta +d;\frac{a\theta +b}{c\theta +d}=\theta ,a,b,c,d\in\Bbb{Z},ad-bc=\pm 1 \}=\{\pm \epsilon_0^n;n\in\Bbb{Z} \}.\]
The number $\epsilon_0$ is  only dependent on $D_\theta$.
\end{fac}
We may assume $\epsilon_0>1$. A real number $\epsilon_0$ is called a fundamental unit. 
It is known that the set 
$\{\pm \epsilon_0^m;m\in\Bbb{Z} \}$ is the unit groups (invertible elements) of $\Bbb{Z}[k\theta]$ 
which is an order of a real quadratic field $\Bbb{Q}(\theta )$. Especially, if either 
$D_\theta\equiv 1$ (mod 4) is square free, or $D_\theta\equiv 8$ or $12$ modulo $16$ and 
$D_\theta /4$ is square free, then $\Bbb{Z}[k\theta]$ is a ring of integers of $\Bbb{Q}(\theta )$. 
We obtain the following proposition.
\begin{pro}
Let $\theta$ be a quadratic number with $k\theta^2+l\theta+m=0$, and let $q$ be a projection in 
$M_n(A_\theta)$ for $n\in\Bbb{N}$. 
Then $qA^n$ is an $A_\theta$-$A_\theta$-equivalence bimodule if and only if $Tr_\theta (q)$ is a
unit (an invertible element) in $\Bbb{Z}[k\theta ]$ which is an order of a real quadratic field 
$\Bbb{Q}(\theta )$.  
\end{pro}
Since $Tr_\theta (q)$ is a positive number, we have
\[\{Tr_\theta (q);qA_\theta^n\;\mathrm{is}\;\mathrm{an}\; A_\theta -A_\theta -\mathrm{equivalence}\;
\mathrm{bimodule} \}=\{\epsilon_0^m;m\in\Bbb{Z}\}.\]
We shall consider the relation between the Picard groups of the irrational rotation algebras 
and the unit groups of the real quadratic fields.
We shall show that $A_\theta$-$A_\theta$-equivalence bimodules are simple form. 
\begin{lem}\label{lem:pro}
Let $\mathcal{E}$ be a $A_\theta$-$A_\theta$-equivalence bimodule. Then there exist
a projection $q\in A_\theta$ and an isomorphism $\phi$ from $A$ to $qA_\theta q$ such 
that $\mathcal{E}$ is isomorphic to $qA_\theta$ with the obvious right $A_\theta$-action, the 
obvious right $A_\theta$-valued inner product, the left $A_\theta$-action such that $a\cdot qb=\phi (a)b$ 
and the left $A_\theta$-valued inner product such that $_{A_\theta}\langle pa,pb\rangle =\phi^{-1}
(pab^*p)$ or $A_\theta q$ with the obvious left $A_\theta$-action, the obvious left $A_\theta$-valued 
inner product, the right $A_\theta$-action such that $pb\cdot a=b\phi (a)$ and the right 
$A_\theta$-valued inner product such that $\langle ap,bp\rangle_{A_\theta} =\phi^{-1}
(pa^*bp)$ for any $a,b\in A_\theta$.
\end{lem}
\begin{proof}
There exist a natural number $n$ and a projection $q\in M_n(A_\theta )$ such that $\mathcal{E}$ 
is isomorphic to $qA_\theta^n$. We denote by $\rho$ an isomorphism from $A_\theta$ to $qM_n(A_\theta )q$  
such that the left $A_\theta$-action of $qA_\theta^n$ is defined by $a\cdot \xi =\rho (a)\xi$ 
for any $a\in A_\theta$ and $\xi\in qA_\theta^n$. \ \\
(1) The case $Tr_\theta (q)\leq 1$.

By Lemma \ref{lem:com}, there exists a unitary element $w$ in $M_n(A_\theta )$ such that 
$w^*qw=q^\prime\otimes e_{11}$ where $e_{11}$ is a rank one projection in $M_n(\Bbb{C})$ and 
$q^\prime$ is a projection in $A_\theta$. Since $w^*qwA_\theta^n$ where 
the left $A_\theta$-action is defined by $a\cdot \xi =w^*\rho (a)w\xi$ for any $a\in A_\theta$ and 
$\xi\in w^*qwA_\theta^n$ is isomorphic to $qA_\theta^n$, $\mathcal{E}$ is isomorphic to 
$(q^\prime\otimes e_{11})A_\theta^n$. It is easy to see that $(q^\prime\otimes e_{11})A_\theta^n$ 
is isomorphic to $q^\prime A_\theta$.\ \\
(2) The case $Tr_\theta (q)\geq 1$.

By Lemma \ref{lem:com}, there exists a unitary element $w$ in $M_n(A_\theta )$ such that 
$w^*qw\geq 1\otimes e_{11}$. It is easy to see that $w^*qwM_n(A_\theta )(1\otimes e_{11})$ 
is isomorphic to $w^*qwA_\theta^n$. Hence $\mathcal{E}$ is isomorphic to 
$w^*qwM_n(A_\theta )(1\otimes e_{11})$. Since $w^*qw(1\otimes e_{11})=1\otimes e_{11}$, 
$\mathcal{E}$ is isomorphic to $w^*qwM_n(A_\theta )w^*qw(1\otimes e_{11})$. Define a map  $\psi$ 
from $A_\theta$ to $w^*qwM_n(A_\theta )w^*qw$ by $\psi (a)=w^*\rho (a)w$ for any $a\in A_\theta$. 
Then $\psi$ is an isomorphism from $A_\theta$ to $w^*qwM_n(A_\theta )w^*qw$. Therefore 
$\mathcal{E}$ is isomorphic to $A_\theta \psi^{-1} (1\otimes e_{11})$ where the right $A_\theta$-action 
is defined by $b\psi^{-1} (1\otimes e_{11})\cdot a=b\psi^{-1} (a\otimes e_{11})$ for any $a,b\in A_\theta$.
\end{proof}
\begin{rem}\label{rem:pro}
(1) By the proof, $A_\theta q_1$ is isomorphic to $q_2A_\theta^n$ such that $Tr_\theta (q_2)=
\tau_\theta (q_1)^{-1}$.\ \\
(2) Lemma \ref{lem:pro} is important in the next section.
\end{rem}
We denote by $\mathcal{E}_{\theta ,\phi ,q}$ (resp. $\mathcal{F}_{\theta ,\phi,q}$) the $A_\theta$-$A_\theta$-equivalence 
bimodule $qA_\theta$ with the left $A_\theta$-action such that $a\cdot qb=\phi (a)b$ 
(resp. $A_\theta q$ with the right $A_\theta$-action such that $bq\cdot a=b\phi (a)$)
for any $a,b\in A_\theta$. We called $\mathcal{E}_{\theta ,\phi ,q}$ and $\mathcal{F}_{\theta ,\phi,q}$ 
real multiplication in the case $q\neq 1$. 
We shall consider the Picard groups of irrational rotation algebras.
\begin{pro}\label{pro:gro}
Let $q_1$ and $q_2$ be projections in $A_\theta$ such that $A_\theta$ is isomorphic to $q_1A_\theta
q_1$ and $q_2A_\theta q_2$. Assume that $\phi_1$ (resp. $\phi_2$) is an isomorphism from $A_\theta$ 
to $q_1A_\theta q_1$ (resp. $q_2A_\theta q_2$). Then $\mathcal{E}_{\theta ,\phi_1,q_1}\otimes\mathcal{E}_
{\theta ,\phi_2,q_2}$ is isomorphic to $\mathcal{E}_{\theta ,\phi_2\circ\phi_1,q}$ where $q$ is a projection in 
$A_\theta$ such that $\tau_\theta (q)=\tau_\theta (q_1)\tau_\theta (q_2)$. Moreover, $\mathcal{E}_{\theta ,\phi_1,q_1}
\otimes\mathcal{F}_{\theta ,\phi_1,q_1}$ and $\mathcal{F}_{\theta ,\phi_1,q_1}\otimes\mathcal{E}_{\theta ,\phi_1,q_1}$ 
are isomorphic to $A_\theta$ with the obvious actions and the inner products.
\end{pro}
\begin{proof}
Let $q=\phi_2 (q_1)$. Then $\tau_\theta (q)=\tau_\theta (q_1)\tau_\theta (q_2)$ and $q\leq q_2$. 
Define a map $F$ from $\mathcal{E}_{\theta ,\phi_1,q_1}\otimes\mathcal{E}_{\theta ,\phi_2,q_2}$ 
to $\mathcal{E}_{\theta ,\phi_2\circ\phi_1,q}$ by $F(q_1a\otimes q_2b)=q\phi_2(a)b$ for any 
$a,b\in A_\theta$ and extend it by the universality. Since $qq_2=q$, if $q_1\otimes q_2a=q_1\otimes q_2b$, 
then $qa=qb$. Hence $F$ is well-defined. Easy computations show that $F$ is an isomorphism of 
equivalence bimodule. Since $\mathcal{F}_{\theta ,\phi_1,q_1}$ is a dual module of 
$\mathcal{E}_{\theta ,\phi_1,q_1}$, $\mathcal{E}_{\theta ,\phi_1,q_1}
\otimes\mathcal{F}_{\theta ,\phi_1,q_1}$ and $\mathcal{F}_{\theta ,\phi_1,q_1}\otimes\mathcal{E}_{\theta ,\phi_1,q_1}$ 
are isomorphic to $A_\theta$ with the obvious actions and the inner products.
\end{proof}
K. Kodaka showed that if $\theta$ is a quadratic number, then the Picard
group of $A_\theta$ is isomorphic to a semidirect product of $\mathrm{Aut} (A_\theta )/
\mathrm{Inn} (A_\theta )$ with $\Bbb{Z}$ \cite{kod2}. By Remark \ref{rem:pro} and Proposition \ref{pro:gro}, 
$\Bbb{Z}$ part is related to the unit group of $\Bbb{Z}[k\theta ]$ which is an order of a real 
quadratic field $\Bbb{Q}(\theta )$. This is shown by another method in \cite{NW}. 

\section{Associated $C^*$-algebras}\label{sec:Ass}
In this section, we construct $C^*$-algebras associated with real multiplication. We show that
these algebras are simple and purely infinite. We recall Cuntz-Pimsner algebras \cite{Pim}. Let $A$ 
be a $C^*$-algebra and $\mathcal{E}$ be a right Hilbert $A$-module. For $\xi ,\eta\in\mathcal{E}$, 
the rank one operator $\Theta_{\xi ,\eta}$ is defined by $\Theta_{\xi ,\eta}(\zeta )=\xi\langle
\eta ,\zeta\rangle_A$ for any $\zeta\in\mathcal{E}$. We denote by $B_A(\mathcal{E})$ 
the algebra of the adjointable operators on $\mathcal{E}$ and by $K_A(\mathcal{E})$ the closure 
of the linear span of rank one operators of $\mathcal{E}$. We say that $\mathcal{E}$ is a Hilbert 
bimodule over $A$ if $\mathcal{E}$ is a right Hilbert $A$-module with a homomorphism $\phi :
A\rightarrow B_A(\mathcal{E})$. We assume that $\mathcal{E}$ is full and $\phi$ is injective. 
We define $I_\mathcal{E}=\phi^{-1}(K_A(\mathcal{E}))$.
The Cuntz-Pimsner algebra $\mathcal{O}_\mathcal{E}$ is the universal $C^*$-algebra generated by 
$A$ and $\{S_\xi ;\xi\in\mathcal{E}\}$ with the following relation 
\[S_{\alpha\xi +\beta\eta}=\alpha S_\xi +\beta S_\eta,\quad aS_\xi b=S_{\phi (a)\xi b},\quad 
S_\xi^*S_\eta =\langle\xi ,\eta\rangle_A\] for any $a,b\in A,\; \xi ,\eta\in\mathcal{E},\; 
\alpha ,\beta\in\Bbb{C}$ and 
\[i_K(\phi (a))=a \] for $a\in I_\mathcal{E}$, where $i_K:I_\mathcal{E}
\rightarrow\mathcal{O}_\mathcal{E}$ is defined by $i_K(\Theta_{\xi ,\eta})=S_\xi S_\eta^*$. 
We shall consider the Cuntz-Pimsner algebras generated by the equivalence bimodules of irrational 
rotation algebras, that are not generated by automorphisms. 
By Lemma \ref{lem:pro}, the equivalence bimodule of irrational rotation algebra
 is isomorphic to $\mathcal{E}_{\theta ,\phi ,q}$ or $\mathcal{F}_{\theta ,\phi ,q}$.
It is easy to see that the Cuntz-Pimsner algebra generated by $\mathcal{E}_{\theta ,\phi ,q}$ is 
isomorphic to the Cuntz-Pimsner algebra generated by $\mathcal{F}_{\theta ,\phi ,q}$. 
We denote by $\mathcal{O}_{qA_\theta ,\phi}$ the Cuntz-Pimsner algebra generated by 
$\mathcal{E}_{\theta ,\phi ,q}$. Since $\mathcal{E}_{\theta ,\phi ,q}$ is an equivalence 
bimodule, $K_{A_\theta}(\mathcal{E}_{\theta ,\phi ,q})$ is isomorphic to $A_\theta$ and 
$S_\xi S_\eta^*=_{A_\theta}\langle\xi ,\eta\rangle$ for $\xi ,\eta\in\mathcal{E}_{\theta ,\phi ,q}$. 
Hence $\mathcal{O}_{qA_\theta, \phi}$ is the universal $C^*$-algebra generated by $A_\theta$ 
and $S_q$ with the following relation; for $a\in A_\theta$, 
\[aS_q=S_q\phi (a),\quad S_q^*S_q=\langle q,q\rangle_{A_\theta}=q,\quad S_qS_q^*
=_{A_\theta}\langle q,q\rangle =\phi^{-1} (q)=1 .\]
Therefore $\mathcal{O}_{qA_\theta \phi}$ is isomorphic to the conner endomorphism crossed 
product $A_\theta\rtimes_\phi\Bbb{N}$. Note that $\phi$ is regarded as an endomorphism such that 
$\phi :A_\theta\rightarrow qA_\theta q\subseteq A_\theta$. 
\begin{thm}
Let $\mathcal{E}_{\theta ,q,\phi}$ be an $A_\theta$-$A_\theta$-equivalence bimodule that is not 
generated by automorphisms. Then the Cuntz-Pimsner algebra $\mathcal{O}_{qA_\theta ,\phi}$ generated 
by $\mathcal{E}_{\theta ,q,\phi}$ is purely infinite, simple, nuclear and in the UCT class. 
\end{thm}
\begin{proof}
By \cite{Ror} (Theorem 3.1.), $\mathcal{O}_{qA_\theta ,\phi}$ is simple and purely infinite because 
$A_\theta$  is a simple unital $C^*$-algebra of real rank zero and with the property of Lemma \ref{lem:com}. 
Since $A_\theta$ is nuclear and in the UCT class, $\mathcal{O}_{qA_\theta ,\phi}$ is nuclear and in the UCT class.
\end{proof}
The isomorphism class of $C^*$-algebra $\mathcal{O}_{qA_\theta ,\phi}$ is completely determined 
by the K-group together with the class of the unit by the classification theorem by Kirchberg-Phillips 
\cite{Kir}, \cite{Phi2}. Especially, the Morita equivalence class of $C^*$-algebra 
$\mathcal{O}_{qA_\theta ,\phi}$ is completely determined by the K-group. 

\section{K-groups}\label{sec:K}
In this section, we compute the K-group of $\mathcal{O}_{qA_\theta ,\phi}$. Let $D_{\mathcal{E}_{\theta ,q,\phi}}$ 
be the linking algebra of $\mathcal{E}_{\theta, q,\phi}$. The linking algebra of $D_{\mathcal{E}_{\theta ,q,\phi}}$ 
is the following form:
\[D_{\mathcal{E}_{\theta ,q,\phi}}=\{\left(
    \begin{array}{cccc}
     x &  \xi \\ 
     \eta &  y
    \end{array}
    \right); x\in K(\mathcal{E}_{\theta ,q,\phi})=qA_\theta q,\;y\in A_\theta ,\;\xi\in\mathcal{E}_{\theta, q,\phi}
,\;\eta\in\mathcal{E}_{\theta, q,\phi}^*\}\]
where $\mathcal{E}_{\theta, q,\phi}^*$ is a dual module of $\mathcal{E}_{\theta, q,\phi}$.
The linking algebra $D_{\mathcal{E}_{\theta ,q,\phi}}$ has the unique normalized trace $Tr_D$ such that 
\[Tr_D(\left(
    \begin{array}{cccc}
     x &  \xi \\ 
     \eta &  y
    \end{array}
    \right))=\frac{\tau_\theta (x) +\tau_\theta (y)}{1+\tau_\theta (q)}.\]
The natural embeddings are denoted by 
\[i_{K(\mathcal{E}_{\theta ,q,\phi})}:K_{A_\theta}(\mathcal{E}_{\theta ,q,\phi})\rightarrow D_{\mathcal{E}_{\theta ,q,\phi}}, \;
i_{A_\theta}:A_\theta\rightarrow D_{\mathcal{E}_{\theta ,q,\phi}},\; i:I_{\mathcal{E}_{\theta ,q,\phi}}
\rightarrow A_\theta .\]
By \cite{kat} (Proposition B.3.), the inclusion $i_{A_\theta}:A_\theta\rightarrow 
D_{\mathcal{E}_{\theta ,q,\phi}}$ induces an isomorphism on the K-groups.
We can define a map $K_*([\mathcal{E}_{\theta ,q,\phi}]):K_*(I_{\mathcal{E}_{\theta ,q,\phi}})\rightarrow 
K_*(A_\theta )$ by the composition of the map $K_*(\phi )$ induced by the restriction of $\phi$ to 
$I_{\mathcal{E}_{\theta ,q,\phi}}$, the map $K_*(i_{K(\mathcal{E}_{\theta ,q,\phi})})$ induced by 
$i_{K(\mathcal{E}_{\theta ,q,\phi})}$ and the inverse of 
the isomorphism $K_*(i_{A_\theta})$. We have the following exact sequence \cite{kat}
\[\begin{CD}
      K_0(I_{\mathcal{E}_{\theta ,q,\phi}})  @>K_0(i)-K_0([\mathcal{E}_{\theta ,q,\phi}])>> K_0(A_\theta ) @>>> K_0(\mathcal{O}_{qA_\theta ,\phi}) \\
               @AAA                         @.                           @VVV                    \\
      K_1(\mathcal{O}_{qA_\theta ,\phi}) @<<<   K_1(A_\theta )  @<<K_1(i)-K_1([\mathcal{E}_{\theta ,q,\phi}])<  K_1(I_{\mathcal{E}_{q,\phi}})
      \end{CD}. \]
Since $\mathcal{E}_{\theta ,q,\phi}$ is an equivalence bimodule, $K_*(i_{K(\mathcal{E}_{\theta ,q,\phi})})$ is an 
isomorphism and $I_{\mathcal{E}_{\theta ,q,\phi}}=A_\theta$.  Therefore $K_*([\mathcal{E}_{\theta ,q,\phi}])\in 
GL(2,\mathbb{Z})$ because $K_0(A_\theta )\cong K_1(A_\theta )\cong\mathbb{Z}^2$ and $K_*(\phi )$ is 
an isomorphism. By \cite{Ell}, there exists an automorphism $\alpha$ of $A_\theta$ such that 
$K_1(\alpha )=g$ for any $g\in GL(2,\mathbb{Z})$. Hence $K_1([\mathcal{E}_{\theta ,q,\phi}])$ is depend on 
$\phi$ and can be any element in $GL(2,\mathbb{Z})$. Let $B_{\theta ,q,\phi}=K_1(i)-K_1([\mathcal{E}_{\theta,q,\phi}])$.
We shall consider $K_0([\mathcal{E}_{\theta ,q,\phi}])$.
\begin{pro}\label{pro:k1}
Let $\mathcal{E}_{\theta ,q,\phi}$ be an $A_\theta$-$A_\theta$-equivalence bimodule such that $\tau_\theta
(q)=c\theta +d$ for $c,d\in\mathbb{Z}$ and $0<\theta <1$. Then there exist integers $a$ and $b$ such that 
$\frac{a\theta +b}{c\theta +d}=\theta$ and 
$K_0([\mathcal{E}_{\theta ,q,\phi}])=\left(
    \begin{array}{cccc}
     d &  b \\ 
     c &  a
    \end{array}
    \right) $ for some basis of $\mathbb{Z}^2$. 
\end{pro}
\begin{proof}
We can choose $<[1],[p]>$ where $\tau_\theta (p)=\theta$ as a basis of $K_0(A_\theta )$ 
\cite{Pim1}, \cite{Rie1}. Since $K_0(i_{A_\theta})$ is an isomorphism, 
$<[i_{A_\theta}(1)],[i_{A_\theta}(p)]>$ is a basis of $K_0(D_{\mathcal{E}_{\theta ,q,\phi}})$.
Easy computations show that $Tr_D (i_{K(\mathcal{E}_{\theta ,q,\phi})}(\phi (1)))=
\frac{\tau_\theta (q)}{1+\tau_\theta (q)}$, $Tr_D (i_{K(\mathcal{E}_{\theta ,q,\phi})}(\phi (p)))
=\frac{\theta\tau_\theta (q)}{1+\tau_\theta (q)}$, $Tr_D (i_{A_\theta}(1))=\frac{1}{1+\tau_\theta (q)}$ 
and $Tr_D(i_{A_\theta}(p))=\frac{\theta}{1+\tau_\theta (q)}$. 
Let $K_0([\mathcal{E}_{\theta ,q,\phi}])=\left(
    \begin{array}{cccc}
     x_{11} &  x_{12} \\ 
     x_{21} &  x_{22}
    \end{array}
    \right)$ for a basis $<[1],[p]>$. 
Then we have 
\[\frac{\tau_\theta (q)}{1+\tau_\theta (q)}=\frac{x_{11}+x_{21}\theta}{1+\tau_\theta (q)},\]
\[\frac{\theta\tau_\theta (q)}{1+\tau_\theta (q)}=\frac{x_{12}+x_{22}\theta}{1+\tau_\theta (q)}.\]
By the discussion in Section \ref{sec:Rea}, there exist integers $a$ and $b$ such that 
$\frac{a\theta +b}{c\theta +d}=\theta$. 
Therefore we have $x_{11}=d$, $x_{12}=c$, $x_{21}=b$ and $x_{22}=a$.
\end{proof}
By the proposition above, we shall show some examples.
\begin{ex}
Let $\theta =\frac{-1+\sqrt{5}}{2}$ and $q$ be a projection in $A_\theta$ such that $\tau_\theta (q)
=\frac{-1+\sqrt{5}}{2}=\theta$.
Then $\frac{-\theta +1}{\theta}=\theta$ and $K_0([\mathcal{E}_{\theta ,q,\phi}])=\left(
    \begin{array}{cccc}
     0 &  1 \\ 
     1 &  -1
    \end{array}
    \right) $. By the exact sequence, we have
\[K_0(\mathcal{O}_{qA_\theta ,\phi})\cong ker(B_{\theta ,q,\phi}),\;
K_1(\mathcal{O}_{qA_\theta ,\phi})\cong \mathbb{Z}^2/Im(B_{\theta ,q,\phi}).\]
\end{ex}
\begin{ex}
Let $\theta =\frac{5+\sqrt{5}}{10}$ and $q$ be a projection in $A_\theta$ such that $\tau_\theta (q)
=\frac{-1+\sqrt{5}}{2}=5\theta -3$.
Then $\frac{2\theta -1}{5\theta -3}=\theta$ and $K_0([\mathcal{E}_{\theta ,q,\phi}])=\left(
    \begin{array}{cccc}
     -3 &  -1 \\ 
     5  &  2
    \end{array}
    \right) $. By the exact sequence, we have
\[K_0(\mathcal{O}_{qA_\theta ,\phi})\cong ker(B_{\theta ,q,\phi}),\;
K_1(\mathcal{O}_{qA_\theta ,\phi})\cong \mathbb{Z}^2/Im(B_{\theta ,q,\phi}).\]
\end{ex}
\begin{ex}\label{ex:1}
Let $\theta =\frac{-1+\sqrt{5}}{2}$ and $q$ be a projection in $A_\theta$ such that $\tau_\theta (q)
=\sqrt{5}-2=2\theta -1$.
Then $\frac{-3\theta +2}{2\theta -1}=\theta$ and $K_0([\mathcal{E}_{\theta ,q,\phi}])=\left(
    \begin{array}{cccc}
     -1 &  2 \\ 
     2 &  -3
    \end{array}
    \right) $. By the exact sequence, we have
\[K_0(\mathcal{O}_{qA_\theta ,\phi})\cong \mathbb{Z}/2\mathbb{Z}\oplus\mathbb{Z}/2\mathbb{Z}
\oplus ker(B_{\theta ,q,\phi}),\;
K_1(\mathcal{O}_{qA_\theta ,\phi})\cong \mathbb{Z}^2/Im(B_{\theta ,q,\phi}).\]
\end{ex}
We shall consider the condition of $\theta_1, \theta_2, q_1, q_2$ that $\mathcal{O}_{q_1A_{\theta_1} ,\phi_1}$ 
is Morita equivalent to $\mathcal{O}_{q_2A_{\theta_2},\phi_2}$. The following proposition is natural.
\begin{pro}\label{pro:Mor}
Let $\mathcal{E}_{\theta_1 ,q_1,\phi_1}$ be an $A_{\theta_1}$-$A_{\theta_1}$-equivalence bimodule and 
$\mathcal{E}_{\theta_2,q_2,\phi_2}$ be an $A_{\theta_2}$-$A_{\theta_2}$-equivalence bimodule such that 
$A_{\theta_1}$ is Morita equivalent to $A_{\theta_2}$ and $\tau_\theta (q_1)=\tau_\theta (q_2)$. 
Assume that $\mathbb{Z}^2/Im((B_{\theta_1,q_1,\phi_1})$ is isomorphic 
to $\mathbb{Z}^2/Im(B_{\theta_2,q_2,\phi_2})$. Then 
$\mathcal{O}_{q_1A_{\theta_1} ,\phi_1}$ is Morita equivalent to $\mathcal{O}_{q_2A_{\theta_2}, \phi_2}$.
\end{pro}
\begin{proof}
We can assume that $0<\theta_1,\theta_2<1$. Since $A_{\theta_1}$ is Morita equivalent to $A_{\theta_2}$, 
there exists $g\in GL(2,\mathbb{Z})$ such that $\theta_1=g\theta_2$ by \cite{Rie1}.  
There exist $c_1,c_2,d_1,d_2\in\mathbb{Z}$ such that $\tau_{\theta_1}(q_1)=\tau_{\theta_2}(q_2)
=c_1\theta_1+d_1=c_2\theta_2+d_2$. The discussion in Section \ref{sec:Rea} shows that there exists 
$a_1,a_2,b_1,b_2\in\mathbb{Z}$ such that $\left(
    \begin{array}{cccc}
     a_1 &  b_1 \\ 
     c_1 &  d_1
    \end{array}
    \right)\theta_1=\theta_1$ and $\left(
    \begin{array}{cccc}
     a_2 &  b_2 \\ 
     c_2 &  d_2
    \end{array}
    \right)\theta_2=\theta_2$.
By the computation, we have $g\left(
    \begin{array}{cccc}
     a_2 &  b_2 \\ 
     c_2 &  d_2
    \end{array}
    \right)g^{-1}=\left(
    \begin{array}{cccc}
     a_1 &  b_1 \\ 
     c_1 &  d_1
    \end{array}
    \right)$.
Hence $Ker(K_0(i)-K_0([\mathcal{E}_{\theta_1 ,q_1,\phi_1}]))$ is isomorphic to $Ker(K_0(i)-K_0
([\mathcal{E}_{\theta_2 ,q_2,\phi_2}]))$ and $Im(K_0(i)-K_0([\mathcal{E}_{\theta_1 ,q_1,\phi_1}]))$ 
is isomorphic to $Im(K_0(i)-K_0([\mathcal{E}_{\theta_2 ,q_2,\phi_2}]))$ by Proposition \ref{pro:k1}. 
Therefore the K-group of $\mathcal{O}_{q_1A_{\theta_1} ,\phi_1}$ is isomorphic to 
the K-group of $\mathcal{O}_{q_2A_{\theta_2} ,\phi_2}$. Consequently, $\mathcal{O}_{q_1A_{\theta_1} ,\phi_1}$ 
is Morita equivalent to $\mathcal{O}_{q_2A_{\theta_2}, \phi_2}$.
\end{proof}
We shall show that the K-group of $\mathcal{O}_{qA_\theta ,\phi}$ is related to Pell's equation.
We review elementary number theory. Let $\theta$ be a quadratic irrational number. 
We may assume that $\theta$ satisfies $k\theta^2+l\theta +m=0$ with a natural number $k$ and 
integers $l,m$ such that $gcd(k,l,m)=1$. The equation is uniquely determined. 
Let $\left(
    \begin{array}{cccc}
     a &  b \\ 
     c &  d
    \end{array}
    \right)$ be an element in $GL(2,\mathbb{Z})$ such that $\frac{a\theta +b}{c\theta +d}=\theta$. 
Then the $\left(
    \begin{array}{cccc}
     a &  b \\ 
     c &  d
    \end{array}
    \right)$ can be written in the form
\[\left(
    \begin{array}{cccc}
     a &  b \\ 
     c &  d
    \end{array}
    \right)=\left(
    \begin{array}{cccc}
     \frac{t+lu}{2} &  ku \\ 
     -mu            &  \frac{t-lu}{2}
    \end{array}
    \right) \]
where $t,u$ are integers such that
\[t^2-D_\theta u^2=4\quad \mathrm{if}\;ad-bc=1 \]
\[t^2-D_\theta u^2=-4\quad \mathrm{if}\;ad-bc=-1. \]
The equation above is called Pell's equation. By the above, $c\theta +d=\frac{t+u\sqrt{D_\theta}}{2}$. 
Especially, if $t,u>0$ are minimum integers satisfying one of the equations above, then 
$\frac{t+u\sqrt{D_\theta}}{2}$ is a fundamental unit.
We shall determine the K-group of $\mathcal{O}_{qA_\theta ,\phi}$.
\begin{thm}
Let $\mathcal{E}_{\theta ,q,\phi}$ be $A_\theta$-$A_\theta$-equivalence bimodule such that 
$\tau_\theta (q)=\frac{t+u\sqrt{D_\theta}}{2}$ where $D_\theta$ is a discriminant of $\theta$. 
Define $B_{\theta ,q,\phi}=K_1(i)-K_1([\mathcal{E}_{\theta,q,\phi}])$.
Then $B_{\theta ,q,\phi}$ can be any element in $\{1-g;g\in GL(2,\mathbb{Z})\}$ 
by the choice of $\phi$ and we have the following.\ \\
(1) If $t^2-D_\theta u^2=4$, then 
\[K_0(\mathcal{O}_{qA_\theta ,\phi})\cong \mathbb{Z}/s_1\mathbb{Z}\oplus
\mathbb{Z}/\frac{2-t}{s_1}\mathbb{Z}\oplus ker(B_{\theta ,q,\phi}),\; 
K_1(\mathcal{O}_{qA_\theta ,\phi})\cong \mathbb{Z}^2/Im(B_{\theta ,q,\phi})\]
where $s_1$ is a maximal integer such that $2-t=es_1^2$ and $u=fs_1$ for some $e,f\in\mathbb{Z}$. \ \\
(2) If $t^2-D_\theta u^2=-4$ then 
\[K_0(\mathcal{O}_{qA_\theta ,\phi})\cong \mathbb{Z}/s_2\mathbb{Z}\oplus
\mathbb{Z}/\frac{t}{s_2}\mathbb{Z}\oplus ker(B_{\theta ,q,\phi}),\; 
K_1(\mathcal{O}_{qA_\theta ,\phi})\cong \mathbb{Z}^2/Im(B_{\theta ,q,\phi})\]
where $s_2$ is a maximal integer such that $t=es_2^2$ and $u=fs_2$ for some $e,f\in\mathbb{Z}$.
\end{thm} 
\begin{proof}
We may consider $0<\theta <1$. Let $\theta^\prime =\left(
    \begin{array}{cccc}
     0 &  1 \\ 
     1 &  1
    \end{array}
    \right)\theta$ and $q^\prime$ be a projection in $A_{\theta^\prime}$ such that $\tau_{\theta^\prime}
(q^\prime )=\tau_\theta (q)$. Then there exists $g\in GL(2,\mathbb{Z})$ such that $K_0(i)-K_0(
[\mathcal{E}_{\theta ,q,\phi}])=g^{-1}(K_0(i)-K_0([\mathcal{E}_{\theta^\prime ,q^\prime ,\phi}]))g$ 
by Proposition \ref{pro:Mor}. 
It is easy to see that $m(\theta^\prime )^2+(l+m)\theta +k+l+m=0$ and $0<\theta^\prime <1$. \ \\
Proof of (1). \ \\
The equation $t^2-D_\theta u^2=4$ gives $e(2+t)=-f^2D_\theta$. It is easy to see that $es_1$ is 
an odd number if and only if $t$ is an odd number and $D_\theta$ is an odd number if and only if $l$ is an 
odd number. From these, we see that $\frac{es_1-lf}{2}$ and $\frac{es_1+lf}{2}$ are integers 
by elementary computations. Therefore there exist 
integers $x_{11},x_{12},x_{21},x_{22}$ such that $K_0(i)-K_0([\mathcal{E}_{\theta ,q,\phi}])=\left(
    \begin{array}{cccc}
     1-\frac{t-lu}{2} &  -ku \\ 
     mu            &  1-\frac{t+lu}{2}
    \end{array}
    \right)=
    s_1\left(
    \begin{array}{cccc}
     x_{11} & x_{12}  \\ 
     x_{21} & x_{22}
    \end{array}
    \right)$.
The greatest common divisor of the matrix elements of $K_0(i)-K_0([\mathcal{E}_{\theta ,q,\phi}])$ 
is equal to the greatest common divisor of the matrix elements of $K_0(i)-K_0([\mathcal{E}_{\theta^\prime ,q^\prime ,\phi}])$ 
because $K_0(i)-K_0([\mathcal{E}_{\theta ,q,\phi}])=
g^{-1}(K_0(i)-K_0([\mathcal{E}_{\theta^\prime ,q^\prime ,\phi}]))g$. 
Since $K_0(i)-K_0([\mathcal{E}_{\theta^\prime ,q^\prime ,\phi}])=\left(
    \begin{array}{cccc}
     1-\frac{t-(l+m)u}{2} &  -mu \\ 
     (k+l+m)u             &  1-\frac{t+(l+m)u}{2}
    \end{array}
    \right)$ and $gcd(k,m,k+l+m)=1$, $gcd(x_{11},x_{12},x_{21},x_{22})=1$.
An easy computation shows that $\mathrm{det}(K_0(i)-K_0([\mathcal{E}_{\theta ,q,\phi}]))=2-t$. 
Hence there exists $g_1,g_2\in GL(2,\mathbb{Z})$ such that 
$g_1(K_0(i)-K_0([\mathcal{E}_{\theta ,q,\phi}]))g_2=\left(
    \begin{array}{cccc}
     s_1 & 0  \\ 
     0 & \frac{2-t}{s_1}
    \end{array}
    \right)$.
By the exact sequence, the proof is complete.\ \\
Proof of (2). \ \\
By an easy computation, $\mathrm{det}(K_0(i)-K_0([\mathcal{E}_{\theta ,q,\phi}]))=-t$. 
It is proved in the similar way of (1). 
\end{proof}
We shall show some examples.
\begin{ex}
Let $\theta =\frac{-1+\sqrt{5}}{2}$ and $q$ be a projection in $A_\theta$ such that $\tau_\theta (q)
=\frac{3-\sqrt{5}}{2}$. Then we have
\[K_0(\mathcal{O}_{qA_\theta ,\phi})\cong ker(B_{\theta ,q,\phi}),\;
K_1(\mathcal{O}_{qA_\theta ,\phi})\cong \mathbb{Z}^2/Im(B_{\theta ,q,\phi}).\]
\end{ex}
\begin{ex}\label{ex:2}
Let $\theta =\sqrt{5}-2$ and $q$ be a projection in $A_\theta$ such that $\tau_\theta (q)
=\sqrt{5}-2=\frac{-4+\sqrt{20}}{2}$. Then we have 
\[K_0(\mathcal{O}_{qA_\theta ,\phi})\cong \mathbb{Z}/4\mathbb{Z}
\oplus ker(B_{\theta ,q,\phi}),\;
K_1(\mathcal{O}_{qA_\theta ,\phi})\cong \mathbb{Z}^2/Im(B_{\theta ,q,\phi}).\]
Therefore $\mathcal{O}_{qA_\theta ,\phi}$ is not Morita equivalent to the $C^*$-algebra 
in Example \ref{ex:1}.
\end{ex}
The following corollary is an extension of Proposition \ref{pro:Mor}.
\begin{cor}\label{cor:mor}
Let $\mathcal{E}_{\theta_1 ,q_1,\phi_1}$ be an $A_{\theta_1}$-$A_{\theta_1}$-equivalence bimodule and 
$\mathcal{E}_{\theta_2,q_2,\phi_2}$ be an $A_{\theta_2}$-$A_{\theta_2}$-equivalence bimodule such that 
$D_{\theta_1}=D_{\theta_2}$ and $\tau_\theta (q_1)=\tau_\theta (q_2)$.
Assume that $\mathbb{Z}^2/Im((B_{\theta_1,q_1,\phi_1})$ is isomorphic to $\mathbb{Z}^2/Im(B_{\theta_2,q_2,\phi_2})$. 
Then $\mathcal{O}_{q_1A_{\theta_1} ,\phi_1}$ is Morita equivalent to $\mathcal{O}_{q_2A_{\theta_2}, \phi_2}$.
\end{cor}
\begin{rem}
Let $\theta_1=\sqrt{10}-3$ and $\theta_2=\frac{2+\sqrt{10}}{3}$. Then $D_{\theta_1}=D_{\theta_2}=40$ 
and $A_{\theta_1}$ is not Morita equivalent to $A_{\theta_2}$.(This fact is related to the fact class 
number of $\mathbb{Q}(\sqrt{10})$ is two.) Therefore Corollary \ref{cor:mor} is an extension 
of Proposition \ref{pro:Mor}. 
\end{rem}

\end{document}